\documentclass[a4paper,12pt]{article}
\usepackage{latexsym,amsmath,amssymb,amsthm}
\theoremstyle{plain}

\newcommand{\rn}{\mathbb{R}^n}

\date{March, 2011}

\title{On Multipliers of Spaces of Harmonic Functions in the
Unit Ball of $\mathbb{R}^n$}
\author{Romi Shamoyan$^1$, Ali Abkar$^2$}
\begin{document}\maketitle
${}^1$Department of Mathematics, Bryansk Technical  University,
Bryansk, Russia, email: rshamoyan@gmail.com\\and

${}^2$Department of Mathematics, Imam Khomeini International
University,
Qazvin 34149, Iran; email: abkar@ikiu.ac.ir\vspace{12pt}\\

 \noindent{\bf Abstract}. We completely describe spaces of multipliers of certain harmonic function spaces
 of Bergman type in $\mathbb{R}^n$. This is the first sharp result of this kind for Bergman type mixed norm spaces of
  harmonic functions in unit ball of $R^n$  \\
 \maketitle\noindent {\bf Key words}: harmonic function spaces,
 multipliers, spherical harmonics.
 \maketitle
\section{Introduction}
Let $B^n$ be the unit ball in Euclidean space $\rn$: that is
$$B^n=\left \{x=(x_1,\dots , x_n)\in \rn :|x|=\left
(\sum_{i=1}^n|x_i|^2\right )^{1/2} <1\right \},$$ and
$S^{n}=\partial B^n$ be the unit sphere in $\rn$; 
$S^{n}=\{ x\in \rn :|x|=1\}$. Let $dm_n(x)$ and $dx^\prime$ be
the normalized Lebesgue measures in $B^n$ and $S^{n}$
respectively. If $x=(x_1,\dots, x_n)\in \rn$, its
projection in $S^{n}$ is defined as usual  by $x=x^\prime |x|$.\\
In [1] the following Bergman spaces $A_\alpha^p(B^n)$ of harmonic
functions in $B^n$ were defined:
\begin{multline*}A_\alpha^p(B^n)=\{f\in h(B^n):\|f\|_{p,\alpha}=\\
\left (\int_0^1\int_{ S^{n}}|f(rx^\prime)|^p(1-r)^\alpha r^{n-1}dr
dx^\prime\right )^{1/p}<\infty \}\end{multline*} and for $p=\infty$
and $0\le \alpha <\infty$ we define
$$A_\alpha^\infty(B^n)=\{f\in h(B^n):\|f\|_{\infty,\alpha}= \sup_{x\in
B^n}|f(x)|(1-|x|)^\alpha<\infty \}$$ where $h(B^n)$ is the space of
all harmonic functions in $B^n$. Note that $A^p_\alpha$ is a Banach
space for $p\ge 1$, and a complete metric space for $0<p< 1$ (see
[1]). In the reference [1], the integral representation was provided
and the existence of bounded projections were proved from Lebesgue
$L^p_\alpha(B^n)$ classes to $A^p_\alpha
(B^n)$ spaces.\\
We will define multipliers of these classes and describe spaces of
multipliers of such type spaces. For this objective, we will need
several additional objects and definitions to formulate the main
result of this paper.\\
It is well-known that every harmonic function $f(x)$ in the unit
ball $B^n$ can be represented as $f(x)=\sum_{k\ge
0}r^kc_ky^{(k)}(x^\prime)$ or
$$f(x)=\sum_{k\ge
0}r^k\left(\sum_{g=1}^{d_k}c_k^{(j)}y_j^{(k)}(x^\prime)\right),\qquad
x=rx^\prime\in B^n,$$ where $y_j^{(k)}(x)$ is a spherical harmonic
function (see [1], [2], [3]). It is well-known that this system of
functions form an orthonormal system on the unit sphere $S^n$ of
$L^2(S^n)$, by taking unions of $y^{(k)}_j(x)$ by $k$, (see [1]). We
will need the following vital kernel function of $2n$ variables $Q_m
(x,y),\, x=rx^\prime , y=\rho y^\prime,\, m>0, \rho >0, r>0$ (see
[1]):
$$Q_m(x,y)=2\sum_{k\ge 0}\frac{\Gamma (m+1+k+n/2}{\Gamma (m +1)\Gamma
(m+n/2)}|x|^k|y|^kZ^{(k)}_{x^\prime}(y^\prime)$$
 where $Z^{(k)}_{x^\prime}(y^\prime)$ is a Zonal harmonic function(see[1]).\\
 \noindent{\bf Definition.} We say that a sequence of complex
 numbers $c_k=\{ c_k^{(j)}: j=1,..., d_k\},\, k\ge 0$ is a
 multiplier from $X$ to $Y$, where $X$ and $Y$ are subspaces of
 $h(B^n)$ if for any $f\in X$ with $$f(x)=\sum_{k\ge
 0}r^k\sum_{j=1}^{d_k}b_k^{(j)}y_j^{(k)}(x^\prime)=
 \sum_{k\ge
 0}r^k\left(b_ky^k(x)\right)$$
 we have the function
 $g(x)=\sum_{k\ge 0}r^kb_kc_ky^k(x\prime)\in Y$ where $x=rx^\prime \in
 B^n$. Indeed
$$g(x)=\sum_{k\ge
 0}r^k\sum_{j=1}^{d_k}b_k^{(j)}c_k^{(j)}y_j^{(k)}(x^\prime).$$
 In this case we will write $\{c_k\}\in M_H(X,Y)$.\\

 We also define general mixed norm classes of harmonic functions in $B^n$ as
 follows: for $0<p,q<\infty ,\, -1<\alpha <\infty$, we set
$$A^{p,q}_\alpha(B^n)=\{f\in h(B^n):\|f\|^p_{p,q,\alpha}<\infty\}$$
where
$$\|f\|^p_{p,q,\alpha}=\int_0^1\left(\int_{S^{n}}|f(|x|x^\prime)|^qdx^\prime \right)^{p/q}(1-|x|^2)^
\alpha |x^{n-1}|d|x|.$$ Note these spaces are Banach spaces for $\min
(p,q)>1$ and are complete metric spaces for $\max (p,q)\le 1$. Note
also that for $0<p<\infty$, we have
$A^{p,p}_\alpha(B^n) =A_\alpha^p(B^n)$.\\
For two harmonic functions $f,g$ with
$$f(x)=\sum_{k\ge 0}r^k\sum_{j=1}^{d_k}c_k^{(j)}y_j^{(k)}(x^\prime);$$ and
$$g(x)=\sum_{k\ge 0}\rho^k\sum_{j=1}^{d_k}b_k^{(j)}y_j^{(k)}(x^\prime)$$
we define the convolution of $f$ and $g$ with
$$(f*g)(x)=\sum_{k\ge 0}r^k\sum_{j=1}^{d_k}c_k^{(j)}b_k^{(j)}y_j^{(k)}(x^\prime);\quad x=rx^\prime.$$
The function
$$P_x(y^\prime)=P(x,y^\prime)=\sum_{k\ge 0}r^kZ_{x^\prime}^{(k)}(y^\prime)=
\sum_{k\ge 0}r^k\left(\sum_{j=1}^{d_k}y_j^{(k)}(y^\prime)y_j^{(k)}(x\prime)\right)$$
is a Poisson kernel(see[1]), and finally
$(\Lambda_{m+1}f)(x)$ is a known fractional derivative of $f(x)$ (see [1]):
$$(\Lambda_{m+1}f)(rx^\prime)=\sum_{k\ge 0}r^kc_ky^{(k)}(x^\prime)
\frac{\Gamma (k+n/2+m+1)}{\Gamma (k+n/2)\Gamma (m
+1)}.$$
In the next section we shall use all these terminologies in more details.
\section{ Main Result}
In this main section we will provide a complete description for the
multipliers of $A_\alpha^{p,1}$ classes with some restrictions on
indexes. This is the first result of this type in harmonic spaces with mixed norm.\\ \noindent{\bf
Theorem.} Let $g(x)$ be a harmonic function in $B^n$ and
$$g(x)=\sum_{k\ge 0}r^k\sum_{j=1}^{d_k}c_k^{(j)}y_j^{(k)}(x^\prime),\quad
x=\rho x^\prime, \rho\in (0,1), x^\prime \in S^{n}.$$ Let $0< p\le
1,\, \alpha\in (0,1),\, m>\max \{\alpha -1, 1/p -1\},\, \beta
>0$. Then the following assertions are equivalent:
 \begin{enumerate} \item [1)]$$ \{c_k^{(j)}: j=1,..., d_k, k\in \mathbb{Z}^+\}
 \in M_H(A_\alpha^{p,1}, A_\beta^{p,1})$$ \item[ 2)]
 $$\sup_{0<\rho<1}\sup_{y^\prime\in S^{n}}\left( \int_{S^{n}}|\Lambda_{m+1}
(g*P_{x^\prime})(\rho y^\prime)| dx^\prime\right)(1-\rho )^{m+1-\alpha +\beta}<\infty.$$
\end{enumerate}
\noindent{\bf Proof.} To prove Theorem, we shall need several rather elementar and partially known  lemmas.\\
\noindent{\bf Lemma 1.} The function $Q_\beta(x,y)$ can be estimated in the following way:
$$|Q_\beta(x,y)|\le \frac{C_1(1-r)^{-\beta}}{|r\rho x^\prime -y^\prime|^{n+[\beta ]}}+
\frac{C_2}{(1-r\rho)^{1+\beta}},$$
where $x=rx^\prime,\, y=\rho y^\prime,\, \beta>-1,\, r,\rho \in (0,1)$. \\
The proof of this lemma can be found in [1].\\
\noindent{\bf Lemma 2.} Let $\alpha>-1$ and $\lambda>\alpha+1$. Then
$$\int_0^1(1-r)^\alpha(1-r\rho)^{-\lambda}dr\le C_{\alpha , \lambda}(1-\rho)^{\alpha -\lambda +1};\qquad
\rho\in(0,1).$$
The proof can be found in [1].\\
\noindent{\bf Lemma 3.} Let $f,g\in h(B^n)$. Then using expansions for f and g mentioned above we have
\begin{multline*}\int_{S^{n}}(g*P_{y})(rx^\prime)f(rx^\prime)dx^\prime=
\int_{S^{n}}\left(\sum_{k\ge 0}r^k\sum_{j=1}^{d_k}C_k^{(j)}y_j^{(k)}(y^\prime)y_j^{(k)}(x^\prime)\right)\\
\left(\sum_{m\ge 0}r^m\sum_{i=1}^{d_m}y_i^{(m)}(x^\prime)b_m^{(i)}dx^\prime\right)=\\
=\sum_{k\ge 0}r^{2k}\sum_{j=1}^{d_k}C_k^{(j)}b_k^{(j)}y_j^{(k)}(y^\prime).
\end{multline*}
The proof is based on the known orthonormality  properties of $y_j^{(k)}$we listed above(see[1]),we omit 
details.\\
\noindent{\bf Lemma 4.} Let $m,k,n\in\mathbb{N}$. Then
$$\int_0^1(1-R^2)^mR^{2k+n-1}dR=\frac{1}{2}\frac{\Gamma(m+1)(\Gamma (k+n/2)}{\Gamma (m+1+n/2 +k)}.$$
The proof can be found in [1].\\
\noindent{\bf Lemma 5.}  Let $0< p\le \infty$,
$$M_p(f,r)=\left(\int_{S^{n}}|f(rx^\prime)|^pdx\prime\right)^{1/p},\qquad 0<r\le 1,$$
where $dx^\prime$ is as we defined above  the Lebesgue measure on $\partial B^n$. Let
$0<q\le 1,\, -1<\beta<\infty$. Then
\begin{multline*}\left( \int_0^1 M_p(f,|y|)\frac{(1-|y|)^\beta}{(1-|x||y|)^{\beta +1}}|y|^{n-1}d|y|\right)^q\\
\le C_q \left( \int_0^1 \frac{M_p^q(f,|y|)(1-|y|)^{\beta
q+q-1}}{(1-|x||y|)^{(\beta +1)q}}|y|^{n-1}d|y|\right).
\end{multline*}
Proof. For each $f\in h(B^n)$, we have $M_p(f,r_1)\le M_p(f,r_2)$
whenever $r_1\le r_2$. Since $0<q\le 1$, it follows that
\begin{multline*} J=
\left( \int_0^1 M_p(f,|y|)\frac{(1-|y|)^\beta}{(1-|x||y|)^{\beta +1}}|y|^{n-1}dy\right)^q\\
\le C_1\sum_{k\ge 0}M_p^q(f, 1-2^{-k-1})(2^{-k\beta
q})\left(\int_{1-2^{-k}}^{1-2^{-k-1}}
\frac{|y|^{n-1}d|y|}{(1-|x||y|)^{\beta+1}}\right)^q\\ \le
C_2\sum_{k\ge 0}M_p^q(f, 1-2^{-k-1})(2^{-k\beta
q})\left(1-|x|(1-2^{-k-1}) \right)^{(\beta+1)q}2^{-kq}.
\end{multline*}
Obviously, $$\frac{1-|x|(1-2^{-k-1})}{1-|x|(1-2^{-k+1})}\ge \frac{1}{4};\quad 0<|x|\le 1,\, k\ge 0.$$
Hence
\begin{multline*}J\le C_3\sum_{k\ge0}\left(\int_{1-2^{-k-1}}^{1-2^{-k-2}}d|y|\right)\left(2^{-k(-1+q+\beta q)}\right)\\
\left(M_p^q(f, 1-2^{-k-1})\right)\left(1-|x|(1-2^{-k+1})\right)^{-(\beta +1)q}\\ \le C_4\sum_{k\ge 0}
\int_{1-2^{-k-1}}^{1-2^{-k-2}}M_p^q(f, |y|)\frac{(1-|y|)^{\beta q+q+1}}{(1-|x||y|)^{(\beta +1)q}}d|y|\\
\le C_5\int_0^1\frac{M_p^q(f,|y|)(1-|y|)^{q(\beta +1)-1}}{(1-|x||y|)^{(\beta+1)q}}d|y|\\ \le
C_6\int_0^1M_p^q(f,t)\frac{(1-t)^{q(\beta +1)-1}}{(1-|x|t)^{(\beta+1)q}}t^{n-1}dt.\end{multline*}
This is what we wanted to prove.\\
\noindent{\bf Lemma 6.} Let $y^\prime\in S^{n}$ be a fixed point in the unit sphere of $\mathbb{R}^n$. Let 
$P(x,y^\prime)=\omega_{n-1}\frac{1-|x|^2}{|x-y^\prime|^n}$ be the Poisson kernel (see[1]) for the unit ball $B^n$, 
$x\in B^n,\, y^\prime\in S^{n}$ and $\omega_{n-1}^{-1}$ is the area of the unit sphere. If $x=rx^\prime$, then
\begin{multline*}P(x,y^\prime)=\sum_{k\ge 0}r^kZ^{(k)}_{x^\prime}(y^\prime)\\
=\sum_{k\ge 0}r^k\left(\sum_{j=1}^{d_k}y_j^{(k)}(y^\prime)y_j^{(k)}(x^\prime)\right)=P_{y^\prime}(rx)\end{multline*}
moreover, for $m\in\mathbb{N}$ and $x=rx^\prime$ we have
\begin{multline*}\int_{S^{n}}(g*P_{y^\prime})(rx^\prime)(f(rx^\prime))dx^\prime=\\
2\int_0^1\int_{S^{n}}\Lambda^{m+1}(g*P_{y^\prime})(rR\xi)(f(rR\xi))(1-R^2)^mR^{n-1}dRd\xi\end{multline*}
where as above $$((\Lambda^{m+1})f)(rx^\prime)=\sum_{k\ge 0}r^kC_k(y^k(x^\prime))
\left(\frac{\Gamma(k+n/2+m+1)}{\Gamma(k+n/2)\Gamma(m+1)}\right).$$
{\bf Proof of Lemma 6.} We have to use the orthonormality of $\{y_j^{(k)}\}_{j=1}^{d_k}$ (see [1]) and Lemma 4.
We have the following chain of equalities
\begin{multline*}2\int_0^1\int_{S^{n}}\Lambda^{m+1}(g*P_{y^\prime})(rR\xi)(f(rR\xi))(1-R^2)^mR^{n-1}dRd\xi\\
=2\int_0^1\int_{S^{n}}\left(\sum_{k\ge 0}(rR)^k\frac{\Gamma(k+n/2+m+1)}{\Gamma(k+n/2)\Gamma(m+1)}
\sum_{j=1}^{d_k}C_k^{(j)}y_j^{(k)}(y^\prime)y_j^{(k)}(\xi)\right)\times \\
\left(\sum_{m\ge 0}(rR)^m\sum_{i=1}^{d_m}b_m^{(i)}y_i^{(m)}(\xi)\right)(1-R^2)^mR^{n-1}dRd\xi\\
=2\int_0^1\sum_{k\ge 0}\frac{\Gamma(k+n/2+m+1)}{\Gamma(k+n/2)\Gamma(m+1)}R^{2k}r^{2k}
\sum_{j=1}^{d_k}C_k^{(j)}b_k^{(j)}y_j^{(k)}(y^\prime)(1-R^2)^mR^{n-1}dR\\
=\sum_{k\ge 0}r^{2k}\sum_{j=1}^{d_k}C_k^{(j)}b_k^{(j)}y_j^{(k)}(y^\prime).
\end{multline*}
On the other hand we have
\begin{multline*}\int_{S^{n}}(g*P_{y^\prime})(rx^\prime)(f(rx^\prime))dx^\prime=\int_{S^{n}}
\left(\sum_{k\ge 0}r^{k}\sum_{j=1}^{d_k}C_k^{(j)}y_j^{(k)}(y^\prime)y_j^{(k)}(x^\prime)\right)\\
\left(\sum_{m\ge 0}r^m\sum_{i=1}^{d_m}y_i^{(m)}(x^\prime)b_m^{(i)}\right
)dx^\prime=\sum_{k\ge 0}r^{2k}\sum_{j=1}^{d_k}C_k^{(j)}b_k^{(j)}y_j^{(k)}(y^\prime).
\end{multline*}
Part two is proved.For part one of lemma see[1]. Now based on these lemmas we are in a position to prove our theorem.\\
{\bf Proof of Theorem.} Fix $y^\prime\in S^{n}$. For each
$y^\prime , m\in\mathbb{N}, m>\alpha -1$ we consider a function
$$Q_m(y,x)=\sum_{k\ge 0}(r^k)\sum_{j=1}^{d_k}\rho^k\left(\frac{\Gamma(k+n/2+m+1)}{\Gamma(k+n/2)
\Gamma(m+1)}y_j^{(k)}(y^\prime)\right ) y_j^{(k)}(x^\prime)$$ where
$x=rx^\prime , y=\rho y^\prime$. Since $\{ C_k^{(j)}: j=1,2,...,d_k,
k\ge 0\}\in M_H(A_\alpha^{p,1}, A_\beta^{p,1})$, we can easily
verify that (see [1])
\begin{multline*}h_y(x)=\sum_{k\ge 0}(r^k)\left(\sum_{j=1}^{d_k}\rho^k\frac{\Gamma(k+n/2+m+1)}{\Gamma(k+n/2)
\Gamma(m+1)}y_j^{(k)}(y^\prime)C_k^{(j)}\right )
y_j^{(k)}(x^\prime)\\
= \sum_{k\ge 0}(r^k\rho
^k)\left(\sum_{j=1}^{d_k}\frac{\Gamma(k+n/2+m+1)}{\Gamma(k+n/2)
\Gamma(m+1)}C_k^{(j)}y_j^{(k)}(x^\prime)\right
)y_j^{(k)}(y^\prime)\end{multline*} is in $A_\beta^{p,1}(B^n)$.
Using the closed graph theorem, for each $m\in\mathbb{N},
m>p(\alpha -1)$ we have
$$\| h_y(x)\|_{A_\beta^{p,1}}\le \| \Lambda
_{m+1}(g*P_{x^\prime})(r\rho y^\prime)\|_{A_\beta^{p,1}}\le C\|
Q_m(y,x)\|_{A_\beta^{p,1}}.$$ Now we use the lemmas proved(Lemma1) to
conclude that
\begin{multline*}
\int_{S^{n}}|Q_\beta (x,y)|dx^\prime \le
C\left(\int_{S^{n}}\frac{(1-r\rho)^{-\beta}}{|r\rho x^\prime
-y^\prime|^{n+[\beta]}}+ \frac{1}{(1-r\rho)^{1+\beta}}\right)\\
\le C \frac{1}{(1-r\rho)^{1+\beta}}.\end{multline*} Hence for
$m>(\alpha -1)p$ we have
\begin{multline*}
\| Q_m(y,x)\|_{A_\alpha^{p,1}}=\left(\int_0^1( \int_{S^{n}}|Q_m
(x,y)|dx^\prime\right) ^p(1-|x|)^{\alpha p-1}|x|^{n-1}dx)^{1/p}\\
\le C(\int_0^1\left(\frac{1}{(1-|x||y|)^{(n+1)p}}(1-|x|)^{\alpha
p-1}|x|^{n-1}d|x|\right)^{1/p}\le C(1-|y|)^{-(m+1-\alpha)}.
\end{multline*}
We now use the fact that for $0<r_1\le r_2<1,\, M_p(f, r_1)\le
M_p(f, r_2)$ to obtain
\begin{multline*}
\int_{S^{n}}| \Lambda _{m+1}(g*P_{x^\prime})(\rho^2
y^\prime)|dx^\prime=\left(\int_{\rho_1}^1
...\int_{\rho_n}^1(1-r)^{\beta p-1}r^{n-1}dr\right)^{-1/p}\\
\int_{\rho_1}^1 ...\int_{\rho_n}^1(1-r)^{\beta p-1}\left(
\int_{S^{n}}| \Lambda _{m+1}(g*P_{x^\prime})(\rho^2
y^\prime)|dx^\prime \cdot r^{n-1}dr\right )^{1/p}\\ \le
C(1-\rho)^{-\beta}\| \Lambda _{m+1}(g*P_{x^\prime})(r\rho
y^\prime)\|_{A^{p,1}_\beta}
\end{multline*}
where $x=rx^\prime ,\, y=\rho y^\prime$. Finally, for $m>\alpha
-1$ we have
\begin{multline*}
\int_{S^{n}}| \Lambda _{m+1}(g*P_{x^\prime})(\rho^2
y^\prime)|dx^\prime\le C(1-\rho)^{-(m+1-\alpha +\beta)}
\end{multline*}
or
\begin{multline*}
(\sup_{0<\rho <1})(\sup_{y^\prime\in S^{n}})\left(
\int_{S^{n}}| \Lambda _{m+1}(g*P_{x^\prime})(\rho^2
y^\prime)|dx^\prime\right) \times (1-\rho)^{m+1-\alpha
+\beta}<\infty. \end{multline*} We now manage to prove the reverse
implication of the theorem. By defintion of convolution above we have
\begin{multline*}
f(x)=\sum_{k\ge
0}r^k\sum_{j=1}^{d_k}c_k^{(j)}y_j^{(k)}(x^\prime),\quad
g(y)=\sum_{m\ge
0}\rho^m\sum_{i=1}^{d_m}b_m^{(i)}y_i^{(m)}(y^\prime),
\\ (f*g)( x)=\sum_{k\ge
0}(r\rho)^k\sum_{j=1}^{d_k}c_k^{(j)}b_k^{(j)}y_j^{(k)}(x^\prime);\quad
x=rx^\prime, y=\rho y^\prime \end{multline*} then it follows from
lemmas we proved
\begin{multline*} h(\rho
rx^\prime)=\int_{S^{n}}(g*P_{x^\prime})(ry^\prime)f(\rho
y^\prime)dy^\prime=\\= 2\int_0^1\int_{S^{n}}(\Lambda
^{m+1}g*P_{x^\prime})(r\rho \xi)(1-R^2)^mR^{n-1}dR d\xi ,\quad
m>\alpha -1.
\end{multline*}
Then we have
\begin{multline*}\int_{S^{n}}|h(\rho rx^\prime )| dx^\prime \le \int_0^1\int_{S^{n}}\int_{S^{n}}
|\Lambda ^{m+1}(g*P_{x^\prime})(rR \xi)|\\
|f(\rho R\xi )|(1-R^2)^mR^{n-1}dR d\xi dx^\prime \le \\ C\int_0^1
\sup_{\xi\in S^{n}}\int_{S^{n}} |\Lambda
^{m+1}(g*P_{x^\prime})(rR \xi)|dx^\prime\\
\times\left(\int_{S^{n}}|f(R\xi )|d\xi\right)(1-R^2)^mR^{n-1}dR.
\end{multline*}
It is well-known that for each fixed $\xi\in S^{n}$ by
subharmonicity, the function $$u(rR\xi, m)= \int_{S^{n}}
|\Lambda^{m+1}(g*P_{x^\prime})(rR \xi)|dx^\prime$$ is growing, hence
$\sup_{\xi\in S^{n}}(u_\xi)$ is also growing ,
and hence we can use the lemmas above (Lemma 5) to get the following
estimate:
\begin{multline*}\int_0^1\left(\int_{S^{n}}|h(rx^\prime )|
dx^\prime\right)^p(1-r)^{\beta p-1}r^{n-1}dr\le \\
C_1\int_0^1\int_0^1 \sup_{\xi\in S^{n}}\left( \int_{S^{n}}
|\Lambda^{m+1}(g*P_{x^\prime})(rR \xi)|dx^\prime\right)^p\\ \times
\left(\int_{S^{n}}|(f(R\xi )|d\xi\right)^p
(1-R)^{pm+p-1}R^{n-1}dR(1-r)^{\beta p-1}r^{n-1}dr.
\end{multline*}
We now use the estimates $(1-R)^{pm}\le(1-rR)^{pm}$ or
$(1-R)^{pm+p-1}\le(1-rR)^{pm+p-1}$ for $m>\frac{1-p}{p}$, and the
conditions imposed in the statement of the theorem to conclude
that
 \begin{multline*}\| h\|^p_{A^{p,1}_\beta} \le
C_2\int_0^1\int_0^1\left( \int_{S^{n}}|f(R\xi)|d\xi \right)^p
(1-Rr)^{p\alpha -p\beta -p}(1-R)^{p-1}\\ R^{(n-1)p}dR(1-r)^{\beta
p-1}r^{n-1}dr \le C_3\int_0^1\left( \int_{S^{n}}|f(R\xi )|d\xi
\right)^p\\
 (1-R)^{p-1}\left(\int_0^1(1-Rr)^{p\alpha -p\beta -p}(1-r)^{\beta
p-1}dr\right )R^{(n-1)p}dR \\ \le  C\| f\|^p_{A_\alpha
^{p,1}},\qquad 0<\alpha<1, \, m\in\mathbb{N};\end{multline*}
 here we used Lemma 2
at the last step. Now the proof of theorem is complete.
We remark finnaly for p=1 this theorem was announced many years ago in [4]

\end{document}